\documentclass[12pt]{article}
\usepackage[dvipsnames]{xcolor}
\usepackage{amsmath}
\usepackage{setspace}
\usepackage{bm}
\usepackage{bbm}
\usepackage{amssymb}
\usepackage{indentfirst}
\usepackage[superscript]{cite}
\usepackage{geometry}
\usepackage{mathtools}
\usepackage{tocloft}

\newcommand*{\citena}[1]{%
\begingroup
[\color{Green}
\romannumeral-`\x 
\setcitestyle{numbers}%
\cite{#1}%
\endgroup
]\ignorespacesafterend
}

\newcommand*{\citesup}[1]{%
\begingroup
\color{Green}
\cite{#1}%
\endgroup
\ignorespacesafterend
}

\newcommand*{\eqrefe}[1]{%
\begingroup
(\color{BrickRed}
\romannumeral-`\x 
\setcitestyle{numbers}%
\ref{eq:#1}%
\endgroup
)\ignorespacesafterend
}

\newcommand*{\secrefe}[1]{%
\begingroup
(\color{Aquamarine}
\romannumeral-`\x 
\setcitestyle{numbers}%
\ref{#1}%
\endgroup
)\ignorespacesafterend
}

\newtagform{Tags}[\textcolor{BrickRed}]{\color{Black}(}{)}

\usepackage[super]{natbib}
\setcitestyle{super}

\newcommand{\ii}{\bm{i}}

\newenvironment{eqleft}
 {\begin{equation}\hspace{0pt}}
 {\hspace{650pt minus 1fil}\end{equation}
\ignorespacesafterend}

\usepackage{colortbl}
\usepackage{tabularx}
\usepackage{makecell}
\definecolor{Gray}{gray}{0.8}
\newcolumntype{g}{>{\columncolor{Gray}}c}

\begin{document}
\title{An Exact Formula for the Prime Counting Function}
\date{May 24, 2019} 
\author{Jose Risomar Sousa}
\maketitle
\usetagform{Tags}

\begin{abstract}
This paper discusses a few main topics in Number Theory, such as the M\"{o}bius function and its generalization, leading up to the derivation of neat power series for the prime counting function, $\pi(x)$, and the prime-power counting function, $J(x)$. Among its main findings, we can cite the extremely useful inversion formula for Dirichlet series (given $F_a(s)$, we know $a(n)$, which implies the Riemann hypothesis, and enabled the creation of a formula for $\pi(x)$ in the first place), and the realization that sums of divisors and the M\"{o}bius function are particular cases of a more general concept. From this result, one concludes that it's not necessary to resort to the zeros of the analytic continuation of the zeta function to obtain $\pi(x)$.
\end{abstract}

\tableofcontents

\section{Introduction}
Many people have devoted time trying to create formulae to generate prime numbers or to count primes numbers, something that at times has bordered on obsession.\\

These formulae don't seem to have a lot of potential to be used as new tools for analyses, so frequently they are mere curiosities or attempts to prove oneself capable of achieving a goal or winning an intellectual challenge. That is even more true if the formula is very complicated, which seems to be the case of most that were discovered to date.\\

But oblivious to the bleak landscape, I was still very curious and eager to find my own, and in the process, I think I may have created a new tool to look at Dirichlet series.\\

In this paper we create the very first power series for the prime counting function, $\pi(x)$, and for the prime-power counting function, $J(x)$, aside from the Riemann prime counting function (which assumes the Riemann hypothesis, which is still officially unproven). Their relative simplicity stems from a property of Dirichlet series, $F_a(s)$, by which their generated function, $a(n)$, can be expressed as the product of an elementary function and a non-elementary function given by a relatively simple power series. \textbf{This property implies that the function generated by $1/\zeta(s)$ is well-determined by the values of $1/\zeta(s)$ at the positive even integers, doesn't depend on $s$ and can be shown to be the M\"{o}bius function, $\mu(n)$, for every $s$ in the series convergence domain, which is what a certain problem equivalent to the Riemann hypothesis requires for its validity.}\\

Though power series are arguably not the most tractable of solutions, they provide insights that may lead to the discovery of better or more useful formulae. Besides, creating even integral representations for these power series is extremely challenging, let alone closed-forms.\\

After I discovered the results that are discussed here, I did some research in the literature and was very surprised about coincidences between things I found and approaches that had been tried by others before me. Namely, concepts such as the M\"{o}bius function, $\mu(n)$, the Von Mangoldt lambda function, $\Lambda(n)$, and so on and so forth. Maybe these are recurring concepts on the study of the patterns of prime numbers.

\section{Indicator function $k$ divides $n$, $\mathbbm{1}_{k|n}$}
In paper \citena{GHN} we introduced the indicator function $k$ divides $n$, noted $\mathbbm{1}_{k|n}$ and defined as 1 if $k$ divides $n$ and 0 otherwise. This function plays a key role throughout this paper. It can be represented by means of elementary functions:
\begin{equation} \label{eq:k_div_n} \nonumber
\mathbbm{1}_{k|n}=\frac{1}{k}\sum_{j=1}^{k}\cos{\frac{2\pi nj}{k}}=\frac{\cos{2\pi n}-1}{2k}+\frac{1}{2k}\sin{2\pi n}\cot{\frac{\pi n}{k}}
\end{equation}\\
\indent However, that closed-form is not very practical to work with. For example, if $k$ divides $n$ we have an undefined product of the type $0\cdot\infty$ (the sine is 0 and the cotangent is $\infty$).\\

\indent Hence, it's more practical to derive a power series for $\mathbbm{1}_{k|n}$, which can be accomplished by expanding the cosine on the left-hand side with Taylor series and employing Faulhaber's formula, as explained in \citena{GHN}, which gives us the below:
\begin{equation} \label{eq:k_div_n_series1} 
\frac{1}{2k}\sin{2\pi n}\cot{\frac{\pi n}{k}}=\sum_{i=0}^{\infty}(-1)^i(2\pi n)^{2i}\sum_{j=0}^{i}\frac{B_{2j} k^{-2j}}{(2j)!(2i+1-2j)!}
\end{equation}

\subsection{Analog of the $\mathbbm{1}_{k|n}$ function}
The analog of the indicator function $\mathbbm{1}_{k|n}$ is the below sum:
\begin{equation} \label{eq:sum_sine1} \nonumber
\frac{1}{k}\sum_{j=1}^{k}\sin{\frac{2\pi nj}{k}}=\frac{\sin{2\pi n}}{2k}+\frac{1}{k}\cot{\frac{\pi n}{k}}\sin^2{\pi n}
\end{equation}\\
\indent As previously, we can obtain a power series for it by expanding the sine with Taylor series and making use of Faulhaber's formula:
\begin{equation} \label{eq:sum_sine4} \nonumber
\frac{1}{2k}\cot{\frac{\pi n}{k}}(1-\cos{2\pi n})=\sum_{i=0}^{\infty}(-1)^i(2\pi n)^{2i+1}\sum_{j=0}^{i}\frac{B_{2j} k^{-2j}}{(2j)!(2i+2-2j)!}
\end{equation}

\subsection{Transforming the power series of $\mathbbm{1}_{k|n}$}
If we remove the first term from the sum in  \eqrefe{k_div_n_series1}, we get the below:
\begin{equation} \label{eq:k_div_n_series2} \nonumber
\mathbbm{1}_{k|n}=\frac{\cos{2\pi n}-1}{2k}+\frac{\sin{2\pi n}}{2\pi n}-
\sum_{i=0}^{\infty}(-1)^i (2\pi n)^{2i+2}\sum_{j=0}^{i} \frac{B_{2j+2}k^{-2j-2}}{(2j+2)!(2i+1-2j)!}
\end{equation}\\
\indent If we keep extracting terms in this fashion, after extracting $q$ terms we would obtain the following equation:
\begin{equation} \label{eq:k_div_n_series4}
\mathbbm{1}_{k|n}=\frac{\cos{2\pi n}-1}{2k}+\frac{\sin{2\pi n}}{2\pi n}\sum_{i=0}^{q-1}\frac{(-1)^i (2\pi n)^{2i} B_{2i}}{(2i)!k^{2i}}+\sum_{i=0}^{\infty} (-1)^{i+q}(2\pi n)^{2i+2q}\sum_{j=0}^{i} \frac{B_{2j+2q}k^{-2j-2q}}{(2j+2q)!(2i+1-2j)!} \nonumber
\end{equation}\\
\indent Since the first two terms are clearly 0 for integer $n$ we can discard them. This means that $\mathbbm{1}_{k|n}$ has an infinite number of alternative power series and can be rewritten as:
\begin{equation} \label{eq:k_div_n_series_gen} 
\mathbbm{1}_{k|n}=\sum_{i=q}^{\infty} (-1)^{i}(2\pi n)^{2i}\sum_{j=q}^{i} \frac{B_{2j}k^{-2j}}{(2j)!(2i+1-2j)!}
\end{equation}

\section{Sum of powers of divisors of $n$}
Let's define the function $\sigma^2_{m}(n)$ as the sum of the $m$-th powers of the integer divisors of $n$ (the  superscript $2$ was chosen for convenience and will make sense when we reach section \secrefe{Dual}). In mathematical notation, for any complex $m$:
\begin{equation} \label{eq:defin} \nonumber
\sigma^2_{m}(n)=\sum_{k|n}k^{m} 
\end{equation}
\indent By using the power series we derived for $\mathbbm{1}_{k|n}$, such as \eqrefe{k_div_n_series1} or \eqrefe{k_div_n_series_gen}, we can obtain a power series for $\sigma^2_{m}(n)$. For reasons that should be apparent soon, form \eqrefe{k_div_n_series1} is preferred.
\subsection{The divisors count function} \label{divsig0}
Let's start by deriving a power series for the number of divisors of an integer $n$, $\sigma^2_{0}(n)$, also known as $d(n)$. If we take equation \eqrefe{k_div_n_series1} and sum $k$ over the positive integers, we get $\sigma^2_{0}(n)$:
\begin{equation} \label{eq:div_sig_0_old}
\sigma^2_{0}(n)=\sum_{k=1}^{\infty}\mathbbm{1}_{k|n}=\sum_{i=0}^{\infty} (-1)^i (2\pi n)^{2i}\sum_{j=0}^{i} \frac{B_{2j}\zeta(2j)}{(2j)!(2i+1-2j)!}
\end{equation}\\
\indent Now, by recalling the closed-form of the zeta function at the even integers:
\begin{equation} \label{eq:zeta_at_even} \nonumber
\zeta(2j)=-\frac{(-1)^j(2\pi)^{2j} B_{2j}}{2(2j)!} \Rightarrow \frac{B_{2j}}{(2j)!}=-2(-1)^j(2\pi)^{-2j}\zeta(2j) \text{,}
\end{equation}\\
\noindent we can replace $B_{2j}/(2j)!$ in equation  \eqrefe{div_sig_0_old} and express $\sigma^2_{0}(n)$ in a more mnemonic form:
\begin{equation} \label{eq:div_sig_0}
\sigma^2_{0}(n)=-2\sum_{i=0}^{\infty} (-1)^{i}(2\pi n)^{2i}\sum_{j=0}^{i}\frac{(-1)^j (2\pi)^{-2j}\zeta(2j)^2}{(2i+1-2j)!}
\end{equation}\\
\indent A similar rationale that also stems from equation \eqrefe{k_div_n_series1} can be applied to obtain $\sigma^2_{m}(n)$ for any complex $m$ (here we leave the equation in its original form):
\begin{equation} \label{eq:div_sig_-2m} \nonumber
\sigma^2_{m}(n)=\sum_{k=1}^{\infty}\mathbbm{1}_{k|n}\cdot k^{m}=\sum_{i=0}^{\infty}(-1)^i (2\pi n)^{2i}\sum_{\substack{j=0 \\ 2j-m\neq 1}}^{i} \frac{B_{2j}\zeta(2j-m)}{(2j)!(2i+1-2j)!}
\end{equation}\\
\indent Note we need to be careful to avoid the zeta function pole, hence $2j-m\neq 1$.\\

And if we multiply $\sigma^2_{m}(n)$ by $n^{-m}$ we obtain $\sigma^2_{-m}(n)$. It's not so obvious to see why this works (if $k$ is a divisor of $n$, $k/n$ is the reciprocal of another divisor of $n$):
\begin{equation} \label{eq:div_sig_+2m} \nonumber
\sigma^2_{-m}(n)=\sum_{k=1}^{\infty}\mathbbm{1}_{k|n}\cdot\left(\frac{k}{n}\right)^{m}=\sum_{i=0}^{\infty} (-1)^i (2\pi)^{2i}n^{2i-m}\sum_{\substack{j=0 \\ 2j-m\neq 1}}^{i}\frac{B_{2j}\zeta(2j-m)}{(2j)!(2i+1-2j)!}
\end{equation}\\
\indent Another way to obtain $\sigma^2_{m}(n)$ can be achieved using \eqrefe{k_div_n_series_gen}:
\begin{equation} \label{eq:div_sig_+2m_a} \nonumber
\sigma^2_{m}(n)=\sum_{k=1}^{\infty}\mathbbm{1}_{k|n}\cdot {k}^{m}=\sum_{k=1}^{\infty}\sum_{i=q}^{\infty}(-1)^{i}(2\pi n)^{2i}\sum_{j=q}^{i}\frac{B_{2j}k^{m-2j}}{(2j)!(2i+1-2j)!} \Rightarrow
\end{equation}
\begin{equation} \label{eq:div_sig_+2m_b} \nonumber
\sigma^2_{m}(n)=\sum_{i=q}^{\infty}(-1)^{i}(2\pi n)^{2i}\sum_{j=q}^{i}\frac{B_{2j}\zeta(2j-m)}{(2j)!(2i+1-2j)!}
\end{equation}\\
\indent And the above holds for positive or negative integer $q$, except that Bernoulli numbers are not defined for negative subscripts. However, they can be analytically continued (and even used to create a formula for $H_{k}(n)$ valid for all complex $k$, as shown in \citena{AC}).\\

As we can see, the function $\mathbbm{1}_{k|n}$ has a lot of interesting properties, and it will be useful on our goal of studying the prime numbers.

\section{Introducing M\"{o}bius $\mu(n)$} \label{Intro_mu}
The formula we created for $\sigma^2_0(n)$ in section \secrefe{divsig0} begs a question: what would happen if we replaced $\zeta(2j)^2$ with $\zeta(2j)^3$? In other words, what does $\sigma^3_0(n)$ give us? 
\begin{equation} \label{eq:div_sig_1} \nonumber
\sigma^3_{0}(n)=-2\sum_{i=0}^{\infty} (-1)^{i}(2\pi n)^{2i}\sum_{j=0}^{i}\frac{(-1)^j (2\pi)^{-2j}\zeta(2j)^3}{(2i+1-2j)!}
\end{equation}\\
\indent To answer this question, we need to rewrite the initial sum that leads to the above formula:
\begin{equation} \label{eq:div_sig_2} \nonumber
\sum_{k_1=1}^{\infty}\sum_{k_2=1}^{\infty}\mathbbm{1}_{k_1\cdot k_2|n}=\sum_{k_1=1}^{\infty}\sum_{k_2=1}^{\infty}\left(\sum_{i=0}^{\infty} (-1)^i (2\pi n)^{2i}\sum_{j=0}^{i} \frac{B_{2j}(k_1\cdot k_2)^{-2j}}{(2j)!(2i+1-2j)!}\right)
\end{equation}\\
\indent Looking at this formula we are led to conclude $\sigma^3_0(n)$ is the number of permutations of elements from the set of divisors of $n$, taken two at a time, that are also divisors of $n$ when multiplied together. It's also equal to the sum of the divisors count function, $d(k)$, over the integers $k$ that divide $n$, as we see in section \secrefe{Dual}.\\

And how about $\sigma^{1}_0(n)$ and $\sigma^{0}_0(n)$? The former is a constant equal to 1 for all integer $n$, and the latter is the sum of the M\"{o}bius function, which we define in the next section, over the integers $k$ that divide $n$. As a preview, $\sigma^{0}_0(n)$ is 1 if $n=1$ and 0 if $n$ is an integer greater than 1, which we prove in section \secrefe{Proof1}.\\

And finally, what is $\sigma^{-1}_0(n)$? That is the M\"{o}bius function itself, $\mu(n)$:
\begin{equation} \label{eq:mu(n)} 
\mu(n)=-2\sum_{i=0}^{\infty} (-1)^{i}(2\pi n)^{2i}\sum_{j=0}^{i}\frac{(-1)^j (2\pi)^{-2j}\zeta(2j)^{-1}}{(2i+1-2j)!} \text{,}
\end{equation}\\
\noindent which we prove in section \secrefe{Proof2}, after we define square-free numbers.

\subsection{Square-free numbers} \label{sq-free}
A square-free number is a number that can't be divided by any squared prime. In other words, if $n$ is square-free, $p_{1}p_{2}\cdots p_{k}$ is its unique prime decomposition. That said, we can define a function $\mu(n)$ such that:\\
\begin{equation}
  \mu(n)=\begin{cases}
    1, & \text{if $n$=1}\\
    (-1)^k, & \text{if $n$ is square-free with $k$ prime factors}\\
    0, & \text{if $n$ is not square-free} \nonumber
  \end{cases}
\end{equation}\\
\indent This function is the M\"{o}bius function from the previous section, which was named after the German mathematician who introduced it.\\

\indent Back to equation \eqrefe{mu(n)}, one of its advantages is the fact that it provides an analytic continuation of $\mu(n)$ onto $\mathbb{C}$.\footnote{Onto is being used loosely here.} It can be rewritten in a different form:
\begin{equation} \label{eq:mu(n)2}
\mu(n)=-\frac{\sin{2\pi n}}{\pi n}\sum_{j=0}^{\infty} \frac{n^{2j}}{\zeta(2j)}
\end{equation}\\
\indent \textbf{The above form has the very same power series expansion as \eqrefe{mu(n)}, but a finite radius of convergence ($|n|<1$), so \eqrefe{mu(n)} is by definition its analytic continuation.}\\

However, this new form is more tractable and useful for performing manipulations in some cases. For example, we can easily find the Taylor series expansion of \eqrefe{mu(n)2} (that is, \eqrefe{mu(n)}), using the Leibniz rule for the $k$-th derivative of the product of two functions:\\
\begin{equation} \nonumber
(f\cdot g)^{(k)}(x)=\sum_{j=0}^{k}{k\choose j}f^{(j)}(x)g^{(k-j)}(x)
\end{equation}\\
\indent Below the graph of $\mu(n)$ was plotted in the $(0,9)$ interval for some insight into its shape and local minima and maxima (it crosses the $x$-axis at the square-full and half-integers):
\begin{center}
\includegraphics{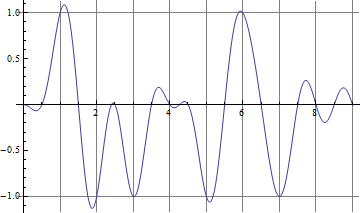}
\end{center}

\indent Note we made the intercept $0$ so the value at $0$ doesn't differ from the other integers (removing successive terms from \eqrefe{mu(n)} or \eqrefe{mu(n)2}, on index $j$, doesn't affect the result for integer or half-integer $n$, as proved in  \eqrefe{k_div_n_series_gen}). It's not easy to apply the Weierstrass factorization theorem here since there are unknown zeros at every level (e.g., $-1$, $0$ or $1$), perhaps even non-real ones. We can remove the half-integer roots dividing $\mu(n)$ by $\cos{\pi n}$ (it can be shown that the resulting power series is \eqrefe{mu(n)} with $2\pi$ replaced by $\pi$), in which case we'd obtain the graph:\\
\begin{center}
\includegraphics{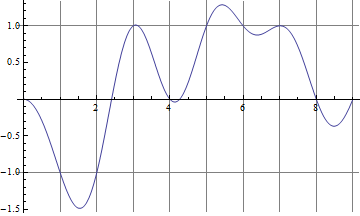}
\end{center}

Going back to equations \eqrefe{mu(n)} and \eqrefe{mu(n)2}, to see an example on why the two forms are the same, we mentioned previously that $\sigma^1_1(n)=1$ for all integer $n$ (in section \secrefe{Dual} we generalize this class of functions). In this case it's possible to produce a closed-form using the generating function of $\zeta(2j)$ that we've created in \citena{GHN}:
\begin{equation} \nonumber
\sigma^1_1(n)=-\frac{\sin{2\pi n}}{\pi n}\sum_{j=0}^{\infty}n^{2j}\zeta(2j)=-\frac{\sin{2\pi n}}{\pi n}\left(-\frac{\pi n\cot{\pi n}}{2}\right)=(\cos{\pi n})^2
\end{equation}\\
\indent The analog of $\sigma^2_{m}(n)$ is $\mu(n)/n^m$, which for complex $m$ is obtained by a simple modification of equation \eqrefe{mu(n)}:
\begin{equation} \nonumber
\frac{\mu(n)}{n^m}=-2\sum_{i=0}^{\infty}(-1)^i (2\pi n)^{2i}\sum_{\substack{j=0 \\ \zeta(2j+m)\neq 0}}^{i}\frac{(-1)^j (2\pi)^{-2j}\zeta(2j+m)^{-1}}{(2i+1-2j)!} \text{,}
\end{equation}
\indent This result is a direct consequence of the inversion formula for Dirichlet series, discussed in section \secrefe{theo}. As before, we need to avoid the zeros of the zeta function if they occur.

\subsection{The Euler product}
The German mathematician Euler discovered an interesting relationship between the zeta function and the primes known as Euler product, which inverted reveals a relationship between the reciprocal of the zeta and the square-free numbers, valid when $\Re{(s)}>1$:
\begin{equation} \label{eq:euler_prod} \nonumber
\frac{1}{\zeta(s)}=\prod_{j=1}^{\infty} \left(1-\frac{1}{p_{j}^{s}}\right)
=\left(1-\frac{1}{2^{s}}\right)\left(1-\frac{1}{3^{s}}\right)\left(1-\frac{1}{5^{s}}\right)\cdots \text{}
\end{equation}\\
\indent Let's denote the set o square-free numbers by $\mathbf{S}$ and let $n$ be a member of $\mathbf{S}$. The inverted Euler product generates all the square-free numbers in $\mathbf{S}$, but each one comes multiplied by $(-1)^k$, where $k$ is the number of primes in the prime decomposition of $n$ ($n=p_{1}p_{2}\cdots p_{k}$).\\

That said, it becomes evident that we can write $\zeta(s)^{-1}$ as a function of $\mu(n)$, if $\Re{(s)}>1$:
\begin{equation} \label{eq:euler_prod2} 
\frac{1}{\zeta(s)}=\sum_{n=1}^{\infty}\frac{\mu(n)}{n^{s}}
\end{equation}

\subsection{Dirichlet series}
\indent The right-hand side of equation \eqrefe{euler_prod2} is one particular example of Dirichlet series. A Dirichlet series is any infinite sum of the type:
\begin{equation} \nonumber
F_a(s)=\sum_{n=1}^{\infty}\frac{a(n)}{n^{s}}
\end{equation}\\
\noindent where $a(n)$ is an arithmetic function and $F_a(s)$ is its generating function.\\

Given two Dirichlet series, $F_a(s)$ and $F_b(s)$, their product is a third Dirichlet series, $F_c(s)$, whose associated function, $c(n)$, is called a Dirichlet convolution of $a$ and $b$, denoted by $c=a*b$.\\

Of particular interest to us, the product of a Dirichlet series $F_a(s)$ and the reciprocal of the zeta function, $\zeta(s)^{-1}$,  is the convolution of $a(n)$ with $\mu(n)$, which is given by the so-called M\"{o}bius inversion formula:
\begin{equation} \label{eq:convo} \nonumber
c(n)=\sum_{k|n}a(k) \mu\left(\frac{n}{k}\right) 
\end{equation}\\
\indent We shall use this result in subsequent proofs.

\subsection{The unit function} \label{Proof1}
Assuming $n$ is integer, let's prove the below assertion, that we referred to in section \secrefe{Intro_mu}:
\begin{equation} \label{eq:test1} \nonumber
\sum_{k|n}\mu(k)=\begin{cases}
      1, & \text{if}\ n=1\\
      0, & \text{if }n>1
\end{cases} 
\end{equation}\\
\indent The proof for the above is pretty simple. If the prime decomposition of $n$ has $k$ primes, then its  number of square-free divisors is $2^k$, half of which have an odd number of prime factors and half of which have an even number of prime factors. Since the former have a negative sign and the latter a positive, they cancel out. The exception is $n=1$, which has no prime factor. $\square$\\

In terms of the function $\mathbbm{1}_{k|n}$ from equation \eqrefe{k_div_n_series1}, or better yet, its alternative form \eqrefe{k_div_n_series_gen}, since equation \eqrefe{euler_prod2} requires $\Re{(s)}>1$, this result implies:
\begin{equation} \nonumber
\sum_{k=1}^{\infty}\mathbbm{1}_{k|n}\cdot\mu(k)=\sum_{i=1}^{\infty} (-1)^i (2\pi n)^{2i}\sum_{j=1}^{i} \frac{B_{2j}}{(2j)!(2i+1-2j)!}\sum_{k=1}^{\infty}\frac{\mu(k)}{k^{2j}} \Rightarrow
\end{equation}
\begin{equation} \nonumber
\mu_0(n)=
\sum_{i=1}^{\infty} (-1)^i (2\pi n)^{2i}\sum_{j=1}^{i} \frac{(-1)^j(2\pi)^{-2j}}{(2i+1-2j)!}=
\begin{cases}
      1, & \text{if}\ n=1\\
      0, & \text{if }n>1
\end{cases}\text{, since \eqrefe{euler_prod2}}\Rightarrow \sum_{k=1}^{\infty}\frac{\mu(k)}{k^{2j}}=\frac{1}{\zeta(2j)}
\end{equation}\\
\indent What this means is that, per the M\"{o}bius inversion formula, the convolution of the function we just called $\mu_0(n)$ with $\mu(n)$ is $\mu(n)$ itself: 
\begin{equation} \nonumber
\mu(n)=\sum_{k|n}\mu_0(k) \mu\left(\frac{n}{k}\right) 
\end{equation}\\
\indent Actually, the convolution of $\mu_0(n)$ with any function is the function itself, and hence $\mu_0(n)$ is known as  the unit function.

\subsection{Cube-free numbers} \label{Proof2}
A cube-free number is a number that can't be divided by any cubed prime. In other words, if $n$ is cube-free, its unique prime decomposition is $n=p_{1}^{e_{1}}p_{2}^{e_{2}}\cdots p_{k}^{e_{k}}$, where the $e_{i}$ are non-zero integer exponents less than or equal to 2. In this context, a prime factor $p_{i}$ is said to be single if $e_{i}=1$, it's said to be double if $e_{i}=2$, and so on. That said, we can introduce a modified M\"{o}bius function of order $2$, $\mu_{2}(n)$, with the following properties:
\begin{equation}
\mu_{2}(n)=\begin{cases} \nonumber
    1, & \text{if $n$=1}\\
    (-2)^{k}, & \text{if $n$ is cube-free with $k$ single prime factors}\\
    0, & \text{if $n$ is not cube-free} 
  \end{cases}
\end{equation}\\
\indent This definition is different from the one proposed by Tom M. Apostol in 1970\citesup{Apostol}, but equal to Popovici's function\citesup{Popovici}.\\
\indent Now, let's take the Euler product from the previous section, and see what it looks like squared:
\begin{equation} \label{eq:euler_prod^2} \nonumber
\frac{1}{\zeta(s)^2}=\prod_{j=1}^{\infty}\left(1-\frac{1}{p_{j}^{s}}\right)^2=\left(1-\frac{1}{2^{s}}\right)^2\left(1-\frac{1}{3^{s}}\right)^2\left(1-\frac{1}{5^{s}}\right)^2\cdots \Rightarrow
\end{equation}

\begin{equation} \label{eq:euler_prod^2_2} \nonumber
\frac{1}{\zeta(s)^2}=\left(1-2\cdot\frac{1}{2^s}+\frac{1}{2^{2s}}\right)\left(1-2\cdot\frac{1}{3^s}+\frac{1}{3^{2s}}\right)\left(1-2\cdot\frac{1}{5^s}+\frac{1}{5^{2s}}\right)\cdots 
\end{equation}\\
\indent Looking at the expansion of $\zeta(s)^{-2}$ above, provided that $\Re{(s)}>1$, it's not very hard to conclude the following equivalence:
\begin{equation} \label{eq:euler_prod3}
\frac{1}{\zeta(s)^2}=\sum_{n=1}^{\infty}\frac{\mu_{2}(n)}{n^{s}}
\end{equation}\\
\indent Now, the convolution of $\mu(n)$ with itself should give us $\mu_2(n)$, after all the latter is generated by $\zeta(s)^{-2}$: 
\begin{equation} \label{eq:theo1_new} \nonumber
\mu_{2}(n)=\sum_{k|n}\mu(k)\mu\left(\frac{n}{k}\right) 
\end{equation}\\
\indent And the above result allows us to state the following theorem:
\begin{eqleft} \label{eq:theo1} \nonumber
\indent \textbf{Theorem 1 } \mu(n)=\sum_{k|n}\mu_{2}(k) 
\end{eqleft}\\
\indent $\textbf{Proof 1}$ This result stems from the M\"{o}bius inversion formula applied to the convolution of $\mu(n)$ with itself.\\

It can also be proved with combinatorics, but we analyze just two possible scenarios. Due to the multiplicative nature of $\mu(n)$ (that is, $\mu(x y)=\mu(x)\mu(y)$, when $x$ and $y$ are co-prime), we can partition $n$ in blocks of co-prime factors where all primes are single, double, and so on, and analyze each one separately. If $n$ has a block of prime factors that are not single, then it's not square-free and hence the sum of $\mu_{2}(k)$ over the integers $k$ that divide $n$ is 0, which we show.\\

For the first scenario, let's assume that $n$ is square-free with $k$ prime factors, $n=p_{1}p_{2}\cdots p_{k}$. Under this scenario we have $2^k$ possible divisors. Let's also assume that $q=p_{1}p_{2}\cdots p_{i}$ is a combination of $i$ out of these $k$ primes. There are $C_{k,i}$ such divisors and they are such that $\mu_{2}(q)=(-2)^i$ if $i>0$. When we plug them into the sum (plus 1, to account for divisor 1) we get:
\begin{equation} 
\sum_{i=0}^{k}C_{k,i}(-2)^i=(1+(-2))^k=(-1)^k=\mu(n) \nonumber
\end{equation}\\
\indent For the second scenario, let's assume that $n$ has $k$ double prime factors, $n=p_{1}^2p_{2}^2\cdots p_{k}^2$. Under this scenario we have $3^k$ possible divisors. Let's also assume that $q=p_{1}p_{2}\cdots p_{i}$, $i>0$, is a combination of $i$ out of these $k$ primes. There are $\sum_{j=0}^{k-i}C_{k,i,j}$ such divisors and they are such that $\mu_{2}(q)=(-2)^i$ if $i>0$, and 1 otherwise.
Let's see why the sum would be 0:
\begin{equation} \nonumber
\sum_{i=0}^{k}\sum_{j=0}^{k-i}C_{k,i,j}(-2)^i =((-2)+1+1)^k=0, \text{ where } C_{k,i,j}=\frac{k!}{i!j!(k-i-j)!} \text{ }\square
\end{equation}\\
\indent Now that theorem 1 has been proved, we can use it to demonstrate the validity of equation \eqrefe{mu(n)}. The demonstration uses equations \eqrefe{k_div_n_series_gen} and \eqrefe{euler_prod3}, and is analogous to the previous one:
\begin{equation} \nonumber
\sum_{k=1}^{\infty}\mathbbm{1}_{k|n}\cdot\mu_2(k)=\sum_{i=1}^{\infty} (-1)^i (2\pi n)^{2i}\sum_{j=1}^{i} \frac{B_{2j}}{(2j)!(2i+1-2j)!}\sum_{k=1}^{\infty}\frac{\mu_2(k)}{k^{2j}} \Rightarrow
\end{equation}
\begin{equation} \nonumber
\mu(n)=\sum_{i=1}^{\infty} (-1)^i (2\pi n)^{2i}\sum_{j=1}^{i} \frac{B_{2j}\zeta(2j)^{-2}}{(2j)!(2i+1-2j)!}=\sum_{i=1}^{\infty} (-1)^i (2\pi n)^{2i}\sum_{j=1}^{i} \frac{(-1)^j(2\pi)^{-2j}\zeta(2j)^{-1}}{(2i+1-2j)!}
\end{equation} 

\subsection{Duality between $\mu_q(n)$ and $\sigma^q_0(n)$} \label{Dual}
From the previous expositions, we can define a generalized M\"{o}bius function of order $q$, $\mu_{q}(n)$, which coincides with Popovici's definition\citesup{Popovici}: $\mu_{q}(n)=\mu *\cdots *\mu$ is the $q$-fold Dirichlet convolution of the M\"{o}bius function with itself. And again, because $\zeta(s)^{-1}$ is the generating function of $\mu(n)$, its convolution with $\mu_q(n)$  justifies the below recurrence:
\begin{equation} \label{eq:tests} \nonumber
\mu_{q+1}(n)=\sum_{k|n}\mu_q(k)\mu\left(\frac{n}{k}\right) \Rightarrow \mu_{q}(n)=\sum_{k|n}\mu_{q+1}(k) \text{, where }\mu_1(n)=\mu(n)\text{.}
\end{equation}\\
\indent Therefore, it follows from this and previous results that:
\begin{equation} \nonumber
\mu_{q}(n)=-2\sum_{i=0}^{\infty} (-1)^i (2\pi n)^{2i}\sum_{j=0}^{i}\frac{(-1)^j (2\pi)^{-2j}\zeta(2j)^{-q}}{(2i+1-2j)!}
\end{equation}\\
\indent The $\mu_q(n)$ expression is insightful, if we think about negative values of $q$. For example:
\begin{equation} \nonumber
\sigma^3_0(n)=\sum_{i=0}^{\infty} (-1)^i (2\pi n)^{2i}\sum_{j=0}^{i} \frac{(-1)^j (2\pi)^{-2j}\zeta(2j)^3}{(2i+1-2j)!}=\sum_{k|n}\sigma^2_0(k)
\end{equation}\\
\indent Again, notice that $\sigma^2_0(k)$ is the number of divisors of $k$, also referred to as $d(k)$.\\

So, we conclude that there is a duality between $\mu_q(n)$ and $\sigma^q_0(n)$, they are equivalent and can be used interchangeably, more precisely, $\mu_{q}(n)=\sigma^{-q}_0(n)$.\\

This finding implies that for positive $q$ the below identities hold, in principle when $\Re{(s)}>1$:
\begin{equation} \nonumber
\frac{1}{\zeta(s)^q}=\sum_{n=1}^{\infty}\frac{\mu_{q}(n)}{n^{s}} \text{, and } \zeta(s)^{q}=\sum_{n=1}^{\infty} \frac{\sigma^q_0(n)}{n^{s}}
\end{equation}\\
\indent The second equation is obvious for case $q=0$ (since $\mu_{0}(n)=\sigma^0_0(n)=0$ for all integer $n$ except 1) and $q=1$ (since $\sigma^1_0(n)=1$ for all integer $n$). 

\section{Inversion formula for Dirichlet series} \label{theo}
We now enunciate a little theorem that relates a Dirichlet series to its coefficients. \\

\indent $\textbf{Theorem 2}$ Suppose that $F_a(s)$ is a Dirichlet series and $a(n)$ is its associated arithmetic function. Then $a(n)$ is given by:
\begin{equation} \label{eq:a(n)} \nonumber
a(n)=-2\sum_{i=0}^{\infty} (-1)^{i}(2\pi n)^{2i}\sum_{j=0}^{i}\frac{(-1)^j (2\pi)^{-2j}F_a(2j)}{(2i+1-2j)!}
\end{equation}\\
\indent $\textbf{Proof 2}$ Although not obvious, this is a very powerful result. The above power series converges for all $n$ and is the analytic continuation of:
\begin{equation} \label{eq:soma}
-\frac{\sin{2\pi n}}{\pi n}\sum_{j=0}^{\infty} n^{2j}F_a(2j) \text{,}
\end{equation}\\
\noindent since they have the same Taylor series expansion and \eqrefe{soma} only converges for $|n|<1$. In some cases it's possible to find a closed-form for $a(n)$, like we did in section \secrefe{sq-free}, though it can be challenging.\\

\indent As for the proof, the $a(n)$ formula is obviously true for $\zeta(s)^q$ for any integer $q$, which we have already proved in the previous sections.\\

\indent For the general case, the proof is sort of trivial: 
\begin{equation} \nonumber
-2\sum_{i=0}^{\infty} (-1)^{i}(2\pi n)^{2i}\sum_{j=0}^{i}\frac{(-1)^j (2\pi)^{-2j}F_a(2j)}{(2i+1-2j)!} \Rightarrow
\end{equation}
\begin{equation} \nonumber
-2\sum_{i=0}^{\infty} (-1)^{i}(2\pi n)^{2i}\sum_{j=0}^{i}\frac{(-1)^j (2\pi)^{-2j}}{(2i+1-2j)!}\sum_{k=1}^{\infty}\frac{a(k)}{k^{2j}} \Rightarrow
\end{equation}
\begin{equation} \nonumber
\sum_{k=1}^{\infty}a(k)\left(-2\sum_{i=0}^{\infty} (-1)^{i}(2\pi n)^{2i}\sum_{j=0}^{i}\frac{(-1)^j (2\pi)^{-2j}k^{-2j}}{(2i+1-2j)!}\right)
\end{equation}\\
\indent The theorem then follows from the following equation:
\begin{equation} \label{eq:delta}
-2\sum_{i=0}^{\infty} (-1)^{i}(2\pi n)^{2i}\sum_{j=0}^{i}\frac{(-1)^j (2\pi k)^{-2j}}{(2i+1-2j)!}=\begin{cases}
      1, & \text{if}\ n=k\\
      0, & \text{otherwise}
\end{cases} 
\end{equation}\\
\indent And the above equation is justified for being the convolution of $\mu_0(n)$ and the associated function of the series $k^{-s}$, $b(n)$, since from the convolution formula:
\begin{equation} \nonumber
c(n)=(\mu_0*b)(n)=\sum_{d|n}\mu_0(d) b\left(\frac{n}{d}\right)=b(n)=\begin{cases}
      1, & \text{if}\ n=k\\
      0, & \text{otherwise}
\end{cases} 
\end{equation}\\
\indent Now, using this result along with the same reasoning employed in the proof from section \secrefe{Proof1}, we derive equation \eqrefe{delta}. $\square$\\

\indent It's quite remarkable that Dirichlet series have coefficients given by Taylor series. One of the advantages of this formula is that if you know $F_a(s)$ at the even integers, you know the coefficients of its series expansion. And if you know the generating function of $F_a(s)$ at the even integers, you know the closed-form of $a(n)$. Another advantage is that it extends $a(n)$ to the complex numbers (for example, primes need no longer be integers). Notice it works even using $F_a(0)$, which is normally a singularity (unless analytically reassigned), but removing this term leaves the result unaltered for integer $n$, as mentioned in section \secrefe{sq-free}.\\

The inversion formula can also be used to check if a function is a Dirichlet series (if at the integers the $a(n)$ are finite and not all zero, the answer is yes).\\

Unfortunately, the inversion formula doesn't apply to the analytic continuation of the Riemann zeta function (as it shouldn't), but it would be really interesting if it did.

\subsection{Riemann hypothesis}
The Riemann hypothesis is equivalent to the statement that the equation:
\begin{equation} \nonumber 
\frac{1}{\zeta(s)}=\sum_{n=1}^{\infty}\frac{\mu(n)}{n^{s}} \text{,}
\end{equation}\\
\noindent is valid for every $s$ with real part greater than $1/2$.\\

\textbf{The inversion formula implies the Riemann hypothesis, since it states that the function generated by $1/\zeta(s)$ is $\mu(n)$ for every $s$ in the series convergence domain.}\\

At this point we have an epiphany, it's quite obvious that the Riemann hypothesis is true. All attempts at claiming to a have proof for it have been and will be met with incredulity and often derision, unless it's a claim by some hotshot. After all, it's been open for more than 150 years and some people refuse to believe that it's that simple. It's quite ironic that we reached a point where proving this conjecture is easier than having it recognized.\\ 

However, the inversion formula settles the riddle, since now we know that the coefficients of a Dirichlet series are the values of a certain analytic function at the positive integers. And in the case of $\zeta(s)^{-1}$, this function happens to be the $\mu(n)$ function exposed before.

\section{Applications}
Even though the possibilities are endless, let's see a few examples.


\subsection{Square root of the zeta}
The function $a(n)$ seems to have a predilection for rational outputs when $F_a(s)$ is some variation of the zeta function:
\begin{equation} \nonumber
\sqrt{\zeta(s)}=1+\frac{1}{2\cdot 2^s}+\frac{1}{2\cdot 3^s}+\frac{3}{8\cdot 4^s}+\frac{1}{2\cdot 5^s}+\frac{1}{4\cdot 6^s}+\frac{1}{2\cdot 7^s} +\frac{5}{16\cdot 8^s}
+\frac{3}{8\cdot 9^s}+\frac{1}{4\cdot 10^s}+\frac{1}{2\cdot 11^s}+\cdots
\end{equation}\\
\indent It's not hard to guess the patterns of $a(n)$: numbers with the same type of prime decomposition have the same coefficients.

\subsection{Zeta raised to $\ii$}
\indent To foray into complex realm, the following coefficients were calculated using the inversion formula:
\begin{equation} \nonumber
\zeta(s)^{\ii}=1+\frac{\ii}{2^s}+\frac{\ii}{3^s}+\frac{-1+\ii}{2\cdot 4^s}+\frac{\ii}{5^s}+\frac{-1}{6^s}+\frac{\ii}{7^s} +\frac{-3+\ii}{6\cdot 8^s}+\frac{-1+\ii}{2\cdot 9^s}+\frac{-1}{10^s}+\frac{\ii}{11^s}+\cdots
\end{equation} 

\subsection{The $n$-th prime number}
If we define the prime zeta function as:
\begin{equation} \nonumber
P(s)=\sum_{n=1}^{\infty}\frac{p_n}{n^s} \text{,}
\end{equation}\\
\noindent then the $n$-th prime is given by:
\begin{equation} \nonumber
p_n=-2\sum_{i=2}^{\infty} (-1)^{i}(2\pi n)^{2i}\sum_{j=2}^{i}\frac{(-1)^j (2\pi)^{-2j}P(2j)}{(2i+1-2j)!}
\end{equation}\\
\indent That's despite our not knowing what the closed-form of $P(s)$ at the even integers is. We need to skip $P(0)$ and $P(2)$ to avoid singularities, since $P(s)$ only converges for $\Re(s)>2$. 

\subsection{Modulus of $\mu(n)$}
\indent Based on the closed-form of the Dirichlet series whose associated function is $|\mu(n)|$, which appears in section \secrefe{sq_free_divs}, we can deduce that: 
\begin{equation} \nonumber
|\mu(n)|=-2\sum_{i=0}^{\infty} (-1)^{i}(2\pi n)^{2i}\sum_{j=0}^{i}\frac{(-1)^j (2\pi)^{-2j}}{(2i+1-2j)!}\frac{\zeta(2j)}{\zeta(4j)}
\end{equation}

\subsection{Liouville function}
If $\lambda(n)=(-1)^{\Omega(n)}$ is the Liouville function, where $\Omega(n)$ is the number of prime factors (with multiplicity) of $n$, then:
\begin{equation} \nonumber
\lambda(n)=-2\sum_{i=0}^{\infty} (-1)^{i}(2\pi n)^{2i}\sum_{j=0}^{i}\frac{(-1)^j (2\pi)^{-2j}}{(2i+1-2j)!}\frac{\zeta(4j)}{\zeta(2j)} \text{, since }\sum_{n=1}^{\infty}\frac{\lambda(n)}{n^{s}}=\frac{\zeta{(2s)}}{\zeta(s)}
\end{equation}

\subsection{Mertens function}
We can extend the Mertens function to the complex domain, based on equation \eqrefe{mu(n)}:
\begin{equation} \nonumber
M(x)=\sum_{n=1}^{x}\mu(n)=-2\sum_{i=0}^{\infty} (-1)^i (2\pi)^{2i}H_{-2i}(x)\sum_{j=0}^{i} \frac{(-1)^j (2\pi)^{-2j}\zeta(2j)^{-1}}{(2i+1-2j)!} \text{,}
\end{equation}\\
\noindent where $H_{-2i}(x)$ is the sum of the $2i$-th powers of the first $x$ positive integers. An asymptotic power series for $M(x)$ is possible using findings from section \secrefe{asym}.

\subsection{Square-free divisors of $n$} \label{sq_free_divs}
Looking back at the topic of integer divisors of $n$, now that square-free numbers and the Euler product have been introduced, we can obtain the sum of powers of square-free divisors of $n$ using the following result from the literature:
\begin{equation} \nonumber
\sum_{n=1}^{\infty}\frac{|\mu(n)|}{n^{s}}=\frac{\zeta{(s)}}{\zeta(2s)} 
\end{equation}\\
\indent Therefore, through the same rationale as before, we conclude that for any complex $m$:
\begin{equation} \nonumber
\sum_{k=1}^{\infty}\mathbbm{1}_{k|n}\cdot |\mu(k)|\cdot k^{m}=\sum_{i=0}^{\infty}(-1)^i (2\pi n)^{2i}\sum_{j=0}^{i}\frac{B_{2j}}{(2j)!(2i+1-2j)!}\frac{\zeta(2j-m)}{\zeta(4j-2m)} \text{, with} \begin{cases} \nonumber 
2j-m\neq 1 & \text{}\\
\zeta(4j-2m) \neq 0 & \text{} 
\end{cases}
\end{equation}\\
\indent In particular, the number of distinct prime factors of $n$ is:
\begin{equation} \nonumber
\log_{2}\left(\sum_{i=0}^{\infty}(-1)^i (2\pi n)^{2i}\sum_{j=0}^{i} \frac{B_{2j}}{(2j)!(2i+1-2j)!}\frac{\zeta(2j)}{\zeta(4j)}\right)
\end{equation}

\subsection{Square root of an integer}
Using theorem 2, we can derive power series expansions for functions that are not analytic at zero (and therefore don't admit a Taylor series at 0), which hold only at the positive integers.\\

For example, if $n$ is a positive integer then:
\begin{equation} \nonumber
\sqrt{n}=-2\sum_{i=0}^{\infty} (-1)^{i}(2\pi n)^{2i}\sum_{j=0}^{i}\frac{(-1)^j (2\pi)^{-2j}}{(2i+1-2j)!}\zeta\left(-\frac{1}{2}+2j\right) \text{, since } \zeta\left(-\frac{1}{2}+s\right)=\sum_{n=1}^{\infty}\frac{\sqrt{n}}{n^{s}}
\end{equation}\\
\indent In fact, the inversion formula posits that for any positive integer $n$ and complex $s$:
\begin{equation} \nonumber
\frac{1}{n^s}=-2\sum_{i=0}^{\infty} (-1)^{i}(2\pi n)^{2i}\sum_{j=0}^{i}\frac{(-1)^j (2\pi)^{-2j}}{(2i+1-2j)!}\zeta\left(s+2j\right) \text{, with }s+2j \neq 1\text{.}
\end{equation}
\indent Hence, it allows us to create a power series for the sum of the square root of the first $n$ positive integers:
\begin{equation} \nonumber
\sum_{k=1}^{n}\sqrt{k}=-2\sum_{i=0}^{\infty} (-1)^i (2\pi)^{2i}H_{-2i}(n)\sum_{j=0}^{i} \frac{(-1)^j (2\pi)^{-2j}}{(2i+1-2j)!}\zeta\left(-\frac{1}{2}+2j\right) \text{}
\end{equation}\\
\indent The downside here is that we don't know the closed-forms for the zeta function at the half-integers.

\subsection{Logarithm of an integer}
Like before, if $n$ and $k$ are positive integers, then:
\begin{equation} \nonumber
(\log{n})^k=-2(-1)^k\sum_{i=0}^{\infty} (-1)^{i}(2\pi n)^{2i}\sum_{j=0}^{i}\frac{(-1)^j (2\pi)^{-2j}\zeta^{(k)}(2j)}{(2i+1-2j)!} \text{, since } \zeta^{(k)}(s)=(-1)^k\sum_{n=1}^{\infty}\frac{(\log{n})^k}{n^{s}}
\end{equation}

\subsection{Von Mangoldt function}
The Von Mangoldt function is defined as:
\begin{equation}
\Lambda(n)=\begin{cases} \nonumber 
\log{p}, & \text{if $n=p^k$ for some prime $p$ and integer k $\geq$ 1}\\
0, & \text{otherwise} 
\end{cases}
\end{equation}\\
\indent We can obtain $\Lambda(n)$ by:
\begin{equation} \nonumber
\Lambda(n)=2\sum_{i=0}^{\infty} (-1)^i (2\pi n)^{2i}\sum_{j=0}^{i} \frac{(-1)^j (2\pi)^{-2j}\zeta'(2j)}{(2i+1-2j)!\zeta(2j)} \text{, since }\frac{\zeta'(s)}{\zeta(s)}=-\sum_{n=1}^{\infty}\frac{\Lambda(n)}{n^{s}}
\end{equation}
\indent From $\Lambda(n)$ we can derive another arithmetic function:
\begin{equation}
\frac{\Lambda(n)}{\log{n}}=\begin{cases} \nonumber 
1/k, & \text{if $n=p^k$ for some prime $p$ and integer k $\geq$ 1}\\
0, & \text{otherwise} 
\end{cases}
\end{equation}\\
\indent Hence, we are able to write $\log{\zeta(s)}$ as a Dirichlet series as follows:
\begin{equation} \nonumber
\frac{1}{\zeta(s)}=\prod_{k=1}^{\infty} \left(1-\frac{1}{p_{k}^{s}}\right) \Rightarrow \log{\zeta(s)}=-\sum_{k=1}^{\infty}\log\left(1-\frac{1}{p_{k}^{s}}\right)=\sum_{k=1}^{\infty}\sum_{i=1}^{\infty}\frac{1}{i \left(p_{k}^i\right)^{s}}=\sum_{n=2}^{\infty}\frac{\Lambda(n)}{\log{n}}\frac{1}{n^{s}}
\end{equation}\\
\indent Therefore, by the inversion formula: 
\begin{equation} \label{eq:van} \nonumber
\frac{\Lambda(n)}{\log{n}}=-2\sum_{i=0}^{\infty} (-1)^i (2\pi n)^{2i}\sum_{j=0}^{i} \frac{(-1)^j (2\pi)^{-2j}\log{\zeta(2j)}}{(2i+1-2j)!} \text{,}
\end{equation}\\
\noindent and this function makes the creation of a formula for $J(x)$ trivial. Plus, when it's combined with $\mu(n)$, we are able to create an exact formula for the prime counting function, $\pi(x)$. Let's rewrite it in its simpler form:
\begin{equation} \label{eq:van2} \nonumber
\frac{\Lambda(n)}{\log{n}}=-\frac{\sin{2\pi n}}{\pi n}\sum_{j=0}^{\infty} n^{2j}\log{\zeta{(2j)}}
\end{equation}

Plotted in the $(0,12)$ interval, we can see that the graph of this function doesn't cross the line $y=1$ only at the primes (note we made the intercept 0 to avoid non-real):
\begin{center}
\includegraphics{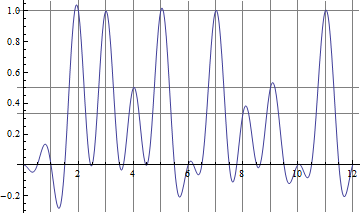}
\end{center}

\indent The maxima of this function are not at the primes, which is why $y=1$ has non-integer roots. It's fair to assume that these roots always occur next to the primes, therefore $\log{\zeta(2)}$ (derived from $y=1$) should be roughly $P(2)$, where $P$ is the usual prime-zeta function.

\newpage

\section{The prime counting function}
We don't even have to know the zeros of the analytic continuation of the zeta function to be able to derive a formula for $\pi(x)$.\\

A number is prime if it's square-free and a prime-power. So, now that we have power series for both $\mu(n)$ and $\Lambda(n)/\log{n}$, we can easily create a function that is 1 whenever $n$ is prime, and 0 otherwise. However, to make the convergence faster, let's divide both by $\cos{\pi n}$ (doing so, we're also removing the half-integer roots):
\begin{equation} \nonumber
\mathbbm{1}_{n\in \mathbb{P}}=-\mu(n)\frac{\Lambda(n)}{\log{n}}/(\cos{\pi n})^2=-\left(-\frac{2\sin{\pi n}}{\pi n}\right)^2\sum_{j=0}^{\infty}\frac{n^{2j}}{\zeta(2j)}\sum_{j=1}^{\infty} n^{2j}\log{\zeta(2j)}
\end{equation}\\
\indent This is not the only and perhaps not even the best way to derive $\mathbbm{1}_{n\in \mathbb{P}}$, but it's the most obvious. If we expand the above function using Taylor series, we arrive at the below series, which has an infinite radius of convergence:
\begin{equation} \nonumber
\mathbbm{1}_{n\in \mathbb{P}}=-\frac{2}{\pi ^2}\left(\sum _{h=1}^{\infty } n^{2h-2}\sum _{j=1}^h \frac{\log\zeta(2j)}{\zeta(2h-2j)}-
\sum _{h=1}^{\infty} n^{2h-2}\sum _{i=1}^h \sum _{j=0}^{h-i} \frac{\log\zeta(2i)}{\zeta(2j)}\frac{(-1)^{h-i-j}(2\pi )^{2h-2i-2j}}{(2h-2i-2j)!}\right)
\end{equation}\\
\indent Fortunately, the above power series can be simplified into a better looking series (which is also more efficient for numeric computation):
\begin{equation} \nonumber
\mathbbm{1}_{n\in \mathbb{P}}=-8\sum _{h=1}^{\infty} n^{2h}\sum _{i=1}^h \log\zeta(2i)\sum _{j=i}^{h}\frac{(-1)^{h-j}(2\pi )^{2h-2j}}{\zeta(2j-2i)(2h+2-2j)!}
\end{equation}\\
\indent To add the half-integer roots back, we just replace $2\pi$ with $4\pi$ in this formula. Finally, the prime counting function is the sum of $\mathbbm{1}_{n\in \mathbb{P}}$ over $n$:
\begin{equation} \nonumber
\pi(x)=-8\sum _{h=1}^{\infty} H_{-2h}(x)\sum _{i=1}^h \log\zeta(2i)\sum _{j=i}^{h}\frac{(-1)^{h-j}(2\pi )^{2h-2j}}{\zeta(2j-2i)(2h+2-2j)!} \text{,}
\end{equation}\\
\noindent where $H_{-2h}(x)$ is the sum of the $2h$-th powers of the first $x$ positive integers, for which the Faulhaber formula provides an analytic continuation.

\newpage 

\subsection{The graph of $\pi(x)$}

Even though it's difficult to compute the  $\pi(x)$ power series for large $x$, the zeros of the zeta function are even harder to compute. Here is a plot of $\pi(x)$ in the $(0,9)$ range (for some reason, the version with half-integer zeros is not as easy to plot):\\
\begin{center}
\includegraphics{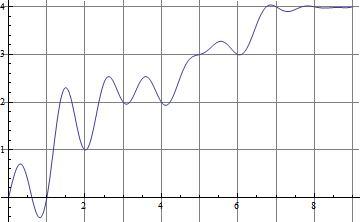}
\end{center}

\subsection{Asymptotic prime counting functions} \label{asym}
If we look back at the formula we just built, it may be possible to simplify it when $x$ is large using the dominant terms of the Faulhaber formula, which is justified by the Euler-Mclaurin formula:
\begin{equation} \nonumber
H_{-2h}(x)=\sum_{n=1}^{x}n^{2h} \sim \frac{x^{2h+1}}{2h+1}+\frac{x^{2h}}{2}
\end{equation}\\
\indent However, given the infinite number of terms of the power series, it's understandable why this approximation should account for only $1/2$ of the total, which is captured in the following conjecture. If $x$ is sufficiently large:
\begin{equation} \nonumber
\pi(x) \sim \frac{\mathbbm{1}_{x\in \mathbb{P}}}{2}-16\sum _{h=1}^{\infty}\frac{x^{2h+1}}{2h+1}\sum _{i=1}^h \log\zeta(2i)\sum _{j=i}^{h}\frac{(-1)^{h-j}(4\pi )^{2h-2j}}{\zeta(2j-2i)(2h+2-2j)!} 
\end{equation}\\
\indent It seems many of the $a(n)$, if not all, have this property. For example, for large $x$, the prime-power counting function, $J(x)$, seems to obey the approximation:
\begin{equation} \nonumber
\sum_{n=2}^{x}\frac{\Lambda(n)}{\log{n}} \sim \frac{\Lambda(x)}{2\log{x}}-4 x\sum_{i=1}^{\infty}\frac{(-1)^i(2\pi x)^{2i}}{2i+1}\sum_{j=1}^{i}\frac{(-1)^j (2\pi)^{-2j}\log{\zeta(2j)}}{(2i+1-2j)!} \text{,}
\end{equation}\\
\noindent for which we've produced the below comparison chart (bigger jumps at the primes):
\begin{table}[h]
\noindent 
\hspace{0cm}
\begin{tabular}{|r|r|r|r|r|r|r|r|r|r|r|r|}

\hline

\multicolumn{1}{|>{\columncolor[gray]{.9}}c|}{$x$} &
\multicolumn{1}{ >{\columncolor[gray]{.9}}c|}{\textbf{Actual}} &
\multicolumn{1}{ >{\columncolor[gray]{.9}}c|}{\textbf{Appr.}} &
\multicolumn{1}{ >{\columncolor[gray]{.9}}c|}{\textbf{\% diff}} &

\multicolumn{1}{ >{\columncolor[gray]{.9}}c|}{$x$} &
\multicolumn{1}{ >{\columncolor[gray]{.9}}c|}{\textbf{Actual}} &
\multicolumn{1}{ >{\columncolor[gray]{.9}}c|}{\textbf{Appr.}} &
\multicolumn{1}{ >{\columncolor[gray]{.9}}c|}{\textbf{\% diff}} &

\multicolumn{1}{ >{\columncolor[gray]{.9}}c|}{$x$} &
\multicolumn{1}{ >{\columncolor[gray]{.9}}c|}{\textbf{Actual}} &
\multicolumn{1}{ >{\columncolor[gray]{.9}}c|}{\textbf{Appr.}} &
\multicolumn{1}{ >{\columncolor[gray]{.9}}c|}{\textbf{\% diff}}  

\\

\hline

1 & 0.00 & 0.00 & 0.0\% & 26 & 11.08 & 11.06 & -0.2\% & 51 & 18.12 & 18.10 & -0.1\% \\ \hline
2 & 1.00 & 1.06 & 5.8\% & 27 & 11.42 & 11.41 & -0.1\% & 52 & 18.12 & 18.14 & 0.1\% \\ \hline
3 & 2.00 & 1.99 & -0.3\% & 28 & 11.42 & 11.43 & 0.1\% & 53 & 19.12 & 19.10 & -0.1\% \\ \hline
4 & 2.50 & 2.48 & -0.7\% & 29 & 12.42 & 12.41 & -0.1\% & 54 & 19.12 & 19.06 & -0.3\% \\ \hline
5 & 3.50 & 3.46 & -1.1\% & 30 & 12.42 & 12.39 & -0.2\% & 55 & 19.12 & 19.09 & -0.2\% \\ \hline
6 & 3.50 & 3.47 & -0.8\% & 31 & 13.42 & 13.38 & -0.3\% & 56 & 19.12 & 19.11 & -0.1\% \\ \hline
7 & 4.50 & 4.48 & -0.4\% & 32 & 13.62 & 13.54 & -0.6\% & 57 & 19.12 & 19.12 & 0.0\% \\ \hline
8 & 4.83 & 4.79 & -1.0\% & 33 & 13.62 & 13.56 & -0.4\% & 58 & 19.12 & 19.16 & 0.2\% \\ \hline
9 & 5.33 & 5.29 & -0.8\% & 34 & 13.62 & 13.59 & -0.2\% & 59 & 20.12 & 20.12 & 0.0\% \\ \hline
10 & 5.33 & 5.32 & -0.2\% & 35 & 13.62 & 13.61 & -0.1\% & 60 & 20.12 & 20.10 & -0.1\% \\ \hline
11 & 6.33 & 6.31 & -0.4\% & 36 & 13.62 & 13.64 & 0.2\% & 61 & 21.12 & 21.08 & -0.2\% \\ \hline
12 & 6.33 & 6.30 & -0.5\% & 37 & 14.62 & 14.61 & -0.1\% & 62 & 21.12 & 21.05 & -0.3\% \\ \hline
13 & 7.33 & 7.29 & -0.6\% & 38 & 14.62 & 14.57 & -0.3\% & 63 & 21.12 & 21.09 & -0.1\% \\ \hline
14 & 7.33 & 7.27 & -0.9\% & 39 & 14.62 & 14.61 & -0.1\% & 64 & 21.28 & 21.27 & -0.1\% \\ \hline
15 & 7.33 & 7.32 & -0.2\% & 40 & 14.62 & 14.65 & 0.2\% & 65 & 21.28 & 21.28 & 0.0\% \\ \hline
16 & 7.58 & 7.60 & 0.3\% & 41 & 15.62 & 15.62 & 0.0\% & 66 & 21.28 & 21.31 & 0.1\% \\ \hline
17 & 8.58 & 8.56 & -0.2\% & 42 & 15.62 & 15.60 & -0.1\% & 67 & 22.28 & 22.27 & 0.0\% \\ \hline
18 & 8.58 & 8.55 & -0.3\% & 43 & 16.62 & 16.58 & -0.2\% & 68 & 22.28 & 22.24 & -0.2\% \\ \hline
19 & 9.58 & 9.54 & -0.4\% & 44 & 16.62 & 16.55 & -0.4\% & 69 & 22.28 & 22.27 & 0.0\% \\ \hline
20 & 9.58 & 9.51 & -0.7\% & 45 & 16.62 & 16.59 & -0.1\% & 70 & 22.28 & 22.31 & 0.1\% \\ \hline
21 & 9.58 & 9.56 & -0.3\% & 46 & 16.62 & 16.63 & 0.1\% & 71 & 23.28 & 23.28 & 0.0\% \\ \hline
22 & 9.58 & 9.60 & 0.2\% & 47 & 17.62 & 17.60 & -0.1\% & 72 & 23.28 & 23.26 & -0.1\% \\ \hline
23 & 10.58 & 10.57 & -0.1\% & 48 & 17.62 & 17.57 & -0.3\% & 73 & 24.28 & 24.24 & -0.2\% \\ \hline
24 & 10.58 & 10.54 & -0.4\% & 49 & 18.12 & 18.08 & -0.2\% & 74 & 24.28 & 24.21 & -0.3\% \\ \hline
25 & 11.08 & 11.05 & -0.3\% & 50 & 18.12 & 18.07 & -0.2\% & 75 & 24.28 & 24.24 & -0.2\% \\ \hline

\end{tabular}
\end{table}

\newpage

Finally, let's have a look at the sum of $d(n)$, and see if it's a just a coincidence. The asymptotic series is:
\begin{equation} \nonumber
\sum_{n=1}^{x}d(n) \sim \frac{d(x)}{2}-4x\sum_{i=1}^{\infty}\frac{(-1)^i(2\pi x)^{2i}}{2i+1}\sum_{j=1}^{i} \frac{(-1)^j (2\pi)^{-2j}\zeta(2j)^{2}}{(2i+1-2j)!} \text{,}
\end{equation}\\
\noindent and a comparison chart is produced below. The conclusion is that this is probably not a coincidence, though the proof might not be simple.

\begin{table}[h]
\noindent 
\hspace{0cm}
\begin{tabular}{|r|r|r|r|r|r|r|r|r|r|r|r|}

\hline

\multicolumn{1}{|>{\columncolor[gray]{.9}}c|}{$x$} &
\multicolumn{1}{ >{\columncolor[gray]{.9}}c|}{\textbf{Actual}} &
\multicolumn{1}{ >{\columncolor[gray]{.9}}c|}{\textbf{Appr.}} &
\multicolumn{1}{ >{\columncolor[gray]{.9}}c|}{\textbf{\% diff}} &

\multicolumn{1}{ >{\columncolor[gray]{.9}}c|}{$x$} &
\multicolumn{1}{ >{\columncolor[gray]{.9}}c|}{\textbf{Actual}} &
\multicolumn{1}{ >{\columncolor[gray]{.9}}c|}{\textbf{Appr.}} &
\multicolumn{1}{ >{\columncolor[gray]{.9}}c|}{\textbf{\% diff}} &

\multicolumn{1}{ >{\columncolor[gray]{.9}}c|}{$x$} &
\multicolumn{1}{ >{\columncolor[gray]{.9}}c|}{\textbf{Actual}} &
\multicolumn{1}{ >{\columncolor[gray]{.9}}c|}{\textbf{Appr.}} &
\multicolumn{1}{ >{\columncolor[gray]{.9}}c|}{\textbf{\% diff}}  

\\

\hline

1 & 1 & 1.16 & 16.5\% & 26 & 91 & 91.26 & 0.3\% & 51 & 211 & 211.22 & 0.1\% \\ \hline
2 & 3 & 3.21 & 6.9\% & 27 & 95 & 95.31 & 0.3\% & 52 & 217 & 217.19 & 0.1\% \\ \hline
3 & 5 & 5.23 & 4.6\% & 28 & 101 & 101.21 & 0.2\% & 53 & 219 & 219.35 & 0.2\% \\ \hline
4 & 8 & 8.22 & 2.7\% & 29 & 103 & 103.31 & 0.3\% & 54 & 227 & 227.36 & 0.2\% \\ \hline
5 & 10 & 10.25 & 2.5\% & 30 & 111 & 111.24 & 0.2\% & 55 & 231 & 231.26 & 0.1\% \\ \hline
6 & 14 & 14.23 & 1.6\% & 31 & 113 & 113.17 & 0.2\% & 56 & 239 & 239.19 & 0.1\% \\ \hline
7 & 16 & 16.23 & 1.4\% & 32 & 119 & 119.28 & 0.2\% & 57 & 243 & 243.06 & 0.0\% \\ \hline
8 & 20 & 20.26 & 1.3\% & 33 & 123 & 123.19 & 0.2\% & 58 & 247 & 247.20 & 0.1\% \\ \hline
9 & 23 & 23.23 & 1.0\% & 34 & 127 & 127.28 & 0.2\% & 59 & 249 & 249.54 & 0.2\% \\ \hline
10 & 27 & 27.21 & 0.8\% & 35 & 131 & 131.42 & 0.3\% & 60 & 261 & 261.26 & 0.1\% \\ \hline
11 & 29 & 29.30 & 1.0\% & 36 & 140 & 140.16 & 0.1\% & 61 & 263 & 262.99 & 0.0\% \\ \hline
12 & 35 & 35.23 & 0.7\% & 37 & 142 & 142.04 & 0.0\% & 62 & 267 & 267.36 & 0.1\% \\ \hline
13 & 37 & 37.18 & 0.5\% & 38 & 146 & 146.31 & 0.2\% & 63 & 273 & 273.39 & 0.1\% \\ \hline
14 & 41 & 41.30 & 0.7\% & 39 & 150 & 150.41 & 0.3\% & 64 & 280 & 280.21 & 0.1\% \\ \hline
15 & 45 & 45.27 & 0.6\% & 40 & 158 & 158.22 & 0.1\% & 65 & 284 & 284.23 & 0.1\% \\ \hline
16 & 50 & 50.18 & 0.4\% & 41 & 160 & 160.23 & 0.1\% & 66 & 292 & 292.13 & 0.0\% \\ \hline
17 & 52 & 52.26 & 0.5\% & 42 & 168 & 168.22 & 0.1\% & 67 & 294 & 294.16 & 0.1\% \\ \hline
18 & 58 & 58.25 & 0.4\% & 43 & 170 & 170.19 & 0.1\% & 68 & 300 & 300.33 & 0.1\% \\ \hline
19 & 60 & 60.26 & 0.4\% & 44 & 176 & 176.35 & 0.2\% & 69 & 304 & 304.36 & 0.1\% \\ \hline
20 & 66 & 66.30 & 0.4\% & 45 & 182 & 182.17 & 0.1\% & 70 & 312 & 312.26 & 0.1\% \\ \hline
21 & 70 & 70.16 & 0.2\% & 46 & 186 & 186.14 & 0.1\% & 71 & 314 & 314.38 & 0.1\% \\ \hline
22 & 74 & 74.19 & 0.3\% & 47 & 188 & 188.46 & 0.2\% & 72 & 326 & 326.19 & 0.1\% \\ \hline
23 & 76 & 76.39 & 0.5\% & 48 & 198 & 198.31 & 0.2\% & 73 & 328 & 327.94 & 0.0\% \\ \hline
24 & 84 & 84.28 & 0.3\% & 49 & 201 & 201.10 & 0.1\% & 74 & 332 & 332.33 & 0.1\% \\ \hline
25 & 87 & 87.11 & 0.1\% & 50 & 207 & 207.23 & 0.1\% & 75 & 338 & 338.34 & 0.1\% \\ \hline

\end{tabular}
\end{table}


\end{document}